\documentstyle[11pt]{article}

\newtheorem{theorem}{Theorem}[section]

\newtheorem{e-proposition}[theorem]{Proposition}
\newtheorem{corollary}[theorem]{Corollary}
\newtheorem{e-definition}[theorem]{Definition\rm}

\begin{document}

 \begin{center}{\bf A Striktpositivstellensatz for measurable functions 
 
 (corrected version)}\bigskip

{Mihai Putinar}\bigskip

{Mathematics Department, 

University of California,

Santa Barbara, CA 93106,

{\rm mputinar@math.ucsb.edu}}
\end{center}

\begin{abstract}
A weighted sums of squares decomposition
of positive Borel measurable functions on a bounded Borel subset
of the Euclidean space is obtained via duality from the spectral
theorem for tuples of commuting self-adjoint operators. The
analogous result for polynomials or certain rational functions was
amply exploited during the last decade in a variety of
applications.

\end{abstract}

\section{Introduction}
To put our main result into the current real algebra context, we
recall below the abstract framework for studying linear
decompositions into weighted sums of squares.

Let $A$ be a commutative algebra with 1, over the rational field.
A {\it quadratic module} $Q \subset A$ is a subset of $A$ such
that $Q+Q \subset Q, 1\in Q$ and $a^2Q\subset A$ for all $a\in A$.
We denote by $Q(F; A)$ or simply $Q(F)$ the quadratic module
generated in $A$ by the set $F$. That is $Q(F; A)$ is the smallest
subset of $A$ which is closed under addition and multiplication by
squares $a^2, \ a \in A$, containing $M$ and the unit $1 \in A$.
If $F$ is finite, we say that the quadratic module is finitely
generated. A quadratic module which is also closed under
multiplication is called a {\it quadratic preordering}. A
quadratic module $Q$ is called archimedean if the constant
function $1$ belong to its algebraic interior, that is, for every
$f \in Q$ there exists $\epsilon >0$ such that $1+tf \in Q$ for
all $0 \leq t \leq \epsilon.$

Assume that $A = {{\mathbf R}}[x_1,...,x_d]$ is the polynomial
algebra. The positivity set $P(Q)$ of $Q \subset {\mathbf
R}[x_1,...,x_d] $ is the set of all points $x \in {\mathbf{R}}^d$
for which $q(x) \geq 0, \  q \in Q$.

The following Striktpositivestellensatz has attracted in the last
decade a lot of attention from practitioners of polynomial
optimization:
{\it Let $Q\subset {{\mathbf R}}[x_1,...,x_d]$ be an
archimedean quadratic module and assume that a
polynomial $f$ is positive on $P(Q)$. Then $f \in Q$.}

This fact was discovered by the author \cite{Put}, generalizing
Schm\"udgen's Striktpositivestellensatz \cite{Sc} for the {\it
finitely generated preordering} associated to a compact
non-negativity set.

In plain language, the above result can be stated as follows.
Denote by $\Sigma^2$ the convex cone of sums of squares in the
polynomials ring ${{\mathbf R}}[x_1,...,x_d]$. Let $p_1,...,p_k \in
{{\mathbf R}}[x_1,...,x_d]$ be polynomials, so that the quadratic
module generated by them $Q(p_1,...,p_k) = \Sigma^2 + p_1 \Sigma^2
+...+ p_k \Sigma^2$ is archimedean. That is, there exists
$\epsilon>0$ so that $ 1-\epsilon (x_1^2+...+x_d^2) \in Q,$ (for
the reduction of $Q$ archimedean to this criterion see \cite{PD}).
The stated Striktpositivstellensatz asserts: if a polynomial $P$
is positive on the set $ S(p_1,...,p_k) =\{ x = (x_1,...,x_d); \
p_i(x) \geq 0, \  1 \leq i \leq k\},$ then $P \in Q(p_1,...,p_k)$.

A simple duality argument implies under the above conditions that
every linear functional $L \in {\mathbf{R}}[x]'$ which is
non-negative on the quadratic module $Q(p_1,...,p_k)$ is
represented by a positive Borel measure $\mu$, supported by the
basic semi-algebraic set $S(p_1,...,p_k)$:
$$ L(f) = \int_{S(p_1,...,p_k)} f d\mu, \ \  f \in
{\mathbf{R}}[x].$$ In general, a quadratic module $Q \subset A$,
where $A$ is an algebra of measurable functions defined on 
${\mathbf R}^d$ has {\it the moment property} if every linear functional
on $A$ which is non-negative on $q$ is represented by a positive measure
supported by $P(Q)$.

The correspondence between the above
Positivestellensatz and the multivariate moment problem with
prescribed compact semi-algebraic supports works fruitfully in
both directions. First, the original proof of the
Positivstellensatz was obtained via the standard moment problem
solution offered by the spectral theorem (see \cite{Put} and
\cite{PD} for an algebraic proof). Second, and more important for
applications, it was J. B. Lasserre \cite{Lasserre0,Lasserre1} who
has interpreted the moments
$$y_\alpha = \int x^\alpha d\mu, \ \  \alpha \in \mathbf{N}^d,$$
of the representing measure as new independent variables and has
obtained a hierarchy of linear, semi-definite optimization problem
converging to the minimization of a given polynomial on a compact
semi-algebraic set. For more details towards applications and
theoretical ramifications see \cite{Bertsimas,HP,posbook,alggeom}.

\section{Main result} The aim of this note is to prove a
natural extension of the polynomial Positivstellensatz to algebras
of Borel measurable functions defined on Euclidean space. Although
a more general statement, on an arbitrary locally compact space or
even on a non-commutative $C^\ast$-algebra is possible to deduce
with similar techniques, we consider that the Euclidean space
setting is: first, the most important for applications, and
second, it contains a specific feature which makes it worth a
separate discussion.

\begin{theorem} Let $Q$ be a countably generated archimedean quadratic module
contained in the algebra \\
${\mathcal{A}} = {\mathbf{R}}[x_1,...,x_d,h_1,...,h_m]$ spanned by the
coordinate functions and by Borel measurable functions
$h_1,...,h_m$ on ${\mathbf{R}}^d$. Assume that $Q$ possesses the moment property. 
If a function $f \in {\mathcal{A}}$
is positive on $P(Q)$, then $f \in Q$.
\end{theorem}

{\bf Proof.} Since $Q$ is archimedean, there exists $\epsilon>0$
such that $1 -\epsilon(x_1^2+...+x_d^2+h_1^2+...+h_m^2) \in Q$.
Thus the positivity set $P(Q)$ is contained in the ball
$x_1^2+...+x_d^2 \leq 1/\epsilon.$ Because $Q$ is countably
generated, the set $P(Q)$ is Borel measurable.

The fact that $Q$ is archimedean as a convex cone, means that for
every $h \in {\mathcal{A}}$ there exists positive constants $c,C$
with the property $C-h, h-c \in Q$. Assume by contradiction that
the function $f$ does not belong to $Q$. According to Marcel Riesz
extension of positive functionals \cite{MRiesz}(known and
rediscovered over the years by many authors
\cite{Kakutani,KR,Ko}), there exists $L \in {\mathcal{A}}'$ so that:
$$ L(f) \leq 0 \leq L(q), \ L(1)>0, \ q \in Q.$$

Next we use Gelfand-Naimark-Segal construction of a Hilbert space
realization of the functional $L$. Specifically, $L(g^2) \geq 0$
for all $g \in {\mathcal{A}}$, and Cauchy-Schwarz' inequality proves
that the set $\mathcal{N}$ of functions $h \in {\mathcal{A}},
L(h)=0,$ is an ideal; whence we can introduce on the quotient
algebra ${\mathcal{A}}/\mathcal{N}$ the positive definite inner
product
$$ \langle g_1, g_2 \rangle = L(g_1 g_2),\ \ g_1,g_2 \in {\mathcal{A}}/\mathcal{N}.$$
Let $\mathcal{H}$ be the Hilbert space completion of
${\mathcal{A}}/\mathcal{N} \otimes_{{\mathbf R}} \mathbf{C}.$ Since
$Q$ is archimedean, the multiplication operators by each generator
$x_1,...,x_d, h_1,...,h_m$ extends by linearity to a tuple of
commuting {\it bounded} self-adjoint operators on $\mathcal{H}$,
denoted by $(X,H) = (X_1,...,X_d,H_1,...,H_m)$, respectively. In
view of the Spectral Theorem \cite{RN}, there exists a positive
measure $\sigma$ on ${\mathbf{R}}^{d+k}$, so that, for all bounded
Borel functions $F(x_1,...,x_d,y_1,...,y_m)$ we have
$$ \langle F(X_1,...,X_d,H_1,...,H_m)1, 1 \rangle= \int_{{\mathbf{R}}^{d+k}} F
d\sigma.$$ From here we deduce that $H_j = h_j(X_1,...,X_d), 1
\leq j \leq m,$ in the sense of Borel functional calculus for
self-adjoint operators \cite{RN}. Hence the measure $\sigma$ is
the push forward of a positive measure on ${\mathbf{R}}^d$ by the
graph map $ x \mapsto (x,h_1(x),...,h_m(x)), \ x \in
{\mathbf{R}}^d$:
$$ \langle F(X_1,...,X_d,H_1,...,H_m)1, 1 \rangle =
\int_{{\mathbf{R}}^{d}} F(x,h_1(x),...,h_m(x))d\mu(x).$$ This shows
that the measure $\mu$ has compact support, contained in the ball
centered at $0$, of radius $1/\epsilon$. Since the algebra $\mathcal A$ was supposed to possess the
moment property, the support of the measure $\mu$ is contained in the positivity set $P(Q)$.

Finally, recall from the statement that $f|_{P(Q)}>0$. On the
other hand, by the construction of the functional we have
$$ \int_{P(Q)} f d\mu = \int f d\mu=  \langle f(X,H) 1, 1 \rangle = L(f) \leq 0.$$
But $\mu({\mathbf{R}}^d)>0$, and thus we reach a contradiction.\\

The reader encounter no complications in specializing the theorem
above and its proof to a finitely generated quadratic module. We
simply state the result.

\begin{corollary}Let $q_1,...,q_n$ be elements of the algebra
 ${{\mathcal{A}}} = {\mathbf{R}}[x_1,...,x_d,h_1,...,h_m]$ generated by
the coordinate functions and by Borel measurable functions
$h_1,...,h_m$ on ${\mathbf{R}}^d$. Let $\Sigma {\mathcal{A}}^2$ denote
the convex cone of sums of squares, and consider the Borel
measurable set
$$ P(q_0, q_1,...,q_n) = \{ x \in {\mathbf{R}}^d; \ q_i(x) \geq 0, \
0
\leq i \leq n\},$$ where $q_0(x) = 1 -
(x_1^2+...+x_d^2+h_1^2+...+h_m^2)$.

If the quadratic module generated by $q_0,...,q_n$ has the moment property, then
a function $f \in {\mathcal{A}}$ is positive on $P(q_0,
q_1,...,q_n)$ satisfies $f \in \Sigma {\mathcal{A}}^2 + q_0 \Sigma
{\mathcal{A}}^2+...+q_n \Sigma {\mathcal{A}}^2.$\\
\end{corollary}

It is very natural to try to extend Lasserre's linearization and relaxation procedure to this
new framework. Specifically that means to consider the mixed moments
$$ y_{\alpha,\beta} = \int x^\alpha h(x)^\beta d\mu, \ \ \alpha
\in \mathbf{N}^d, \ \beta \in \mathbf{N}^m,$$ as new variables, together will all {\it algebraic} dependence relations
among the
functions \\
$x_1,...,x_d,
h_1(x),...,h_m(x)$. For instance it may
happen that $h_i(x)$ is a characteristic function of a Borel set,
in which case $h_i^2 = h_i$, or that $h(x)$ is an $m$-tuple of
algebraic functions, in which case a polynomial dependence
$P(x,h(x)) = 0$ holds. When no algebraic dependence exists
between the generators $x_1,...,x_d,
h_1(x),...,h_m(x),$ this relaxation method will treat them as independent variables.

More details and examples can be found in the recent preprint \cite{LP}.

\end{document}